\documentstyle{amsppt}
\magnification=\magstep1
\pagewidth{4.5in}
\pageheight{7.0in}
\hyphenation{co-deter-min-ant co-deter-min-ants pa-ra-met-rised
pre-print pro-pa-gat-ing pro-pa-gate
fel-low-ship Cox-et-er dis-trib-ut-ive}
\def\leaderfill{\leaders\hbox to 1em{\hss.\hss}\hfill}
\def\A{{\Cal A}}

\def\H{{\Cal H}}
\def\L{{\Cal L}}

\def\e{{\varepsilon}}

\def\P{{\widetilde P}}
\def\Q{{\widetilde Q}}

\def\b0{\text{\bf 0}}

\def\ra{{\ \longrightarrow \ }}

\def\lan{{\langle}}
\def\ran{{\rangle}}

\def\complex{{\Bbb C}}
\def\zed{{\Bbb Z}}

\def\enn{{\Bbb N}}

\def\boxit#1{\vbox{\hrule\hbox{\vrule \kern3pt
\vbox{\kern3pt\hbox{#1}\kern3pt}\kern3pt\vrule}\hrule}}
\def\rabbit{\vbox{\hbox{\kern0pt
\vbox{\kern0pt{\hbox{---}}\kern3.5pt}}}}

\def\tableau#1{
        \hbox {
                \hskip -10pt plus0pt minus0pt
                \raise\baselineskip\hbox{
                \offinterlineskip
                \hbox{#1}}
                \hskip0.25em
        }
}

\def\tabCol#1{
\hbox{\vtop{\hrule
\halign{\strut\vrule\hskip0.5em##\hskip0.5em\hfill\vrule\cr\lower0pt
\hbox\bgroup$#1$\egroup \cr}
\hrule
} } \hskip -10.5pt plus0pt minus0pt}

\def\CR{
        $\egroup\cr
        \noalign{\hrule}
        \lower0pt\hbox\bgroup$
}



\def\blank#1#2{
\hbox to #1{\hfill \vbox to #2{\vfill}}
}

\topmatter
\title Canonical Bases for Hecke Algebra Quotients
\endtitle

\author R.M. Green and J. Losonczy \endauthor
\affil 
{\it Dedicated to Professor R.W. Carter on the occasion of his 65th
birthday}\\
\newline
Department of Mathematics and Statistics\\ Lancaster University\\
Lancaster LA1 4YF\\ England\\
{\it  E-mail:} r.m.green\@lancaster.ac.uk\\
\newline
Department of Mathematics\\
Long Island University\\
Brookville, NY  11548\\
USA\\
{\it  E-mail:} losonczy\@e-math.ams.org\\
\endaffil

\abstract
We establish the existence of an IC basis for
the generalized Temperley--Lieb algebra associated to a Coxeter system
of arbitrary type.  We determine this basis explicitly in the case
where the Coxeter system is simply laced and the algebra is finite
dimensional.
\endabstract

\thanks
The first author was supported in part by an award from the Nuffield
Foundation.
\endthanks

\endtopmatter


\centerline{\bf To appear in Mathematical Research Letters}

\head Introduction \endhead
An important construction in the theory of quantum groups and quantum
algebras is that of canonical bases.  
The original example of this construction is the
Kazhdan--Lusztig basis of the Hecke algebra associated to a Coxeter
system, which first appeared in \cite{{\bf 11}}.  
Another well known example
is the canonical basis for $U^+$, the ``plus part'' of the
quantized enveloping algebra associated to a semisimple Lie algebra
over $\complex$; this was discovered independently by Kashiwara
\cite{{\bf 10}} and Lusztig \cite{{\bf 12}}.  
In each of these
examples, the basis which arises has many deep and beautiful
properties, some of which have geometric interpretations.

The general theory of such bases is defined for an $\A$-module (where
$\A$ is the ring of Laurent polynomials $\zed[v, v^{-1}]$) equipped
with an involutive $\zed$-linear map that sends $v$ to $v^{-1}$.  This
theory was developed by Du in \cite{{\bf 2}}, where the bases arising are 
called IC bases.  The
letters ``IC'' stand for ``intersection cohomology''; the name
alludes to the fact that many of the natural examples have interpretations in
terms of perverse sheaves.  However, the existence and uniqueness 
of such a basis is not guaranteed.

The Temperley--Lieb algebra, a finite dimensional algebra
arising in statistical mechanics \cite{{\bf 13}} and knot theory
\cite{{\bf 9}}, may be defined as a certain quotient of the Hecke algebra of a
Coxeter system of type $A$.  It is possible to generalize this
construction to an arbitrary Coxeter system, obtaining the so-called
``generalized Temperley--Lieb algebras'' as quotients.  In this paper,
we show that IC bases exist for these Hecke algebra quotients
associated to a Coxeter system of arbitrary type.
We also determine the bases explicitly
in the case of Coxeter systems of types $A$, $D$
and $E$; it turns out that in each of these cases the basis coincides
with a previously familiar
basis for the corresponding generalized Temperley--Lieb algebra.  The
situation for non-simply-laced Coxeter systems turns out to be 
more complicated and surprising, as we will discuss.

It is tempting to think that IC bases for generalized Temperley--Lieb
algebras may be obtained from the well-known Kazhdan--Lusztig bases
for the Hecke algebra by projection to the quotient, but this is in
fact far from clear.  In general, the kernel of the canonical map 
from the Hecke algebra to the generalized Temperley--Lieb algebra is 
not spanned by the Kazhdan--Lusztig basis
elements which it contains, even for Coxeter systems of low rank such as
type $D_4$, and this causes complications.
However, the situation in type $A$ is relatively simple, as
C.K. Fan and the first author have shown \cite{{\bf 4}}.

Our results give rise to some interesting problems concerning the general
properties of IC bases for generalized Temperley--Lieb algebras; we 
mention a few of these in our remarks.

\vfill\eject

\head 1. Generalized Temperley--Lieb algebras \endhead

Let $X$ be a Coxeter graph, of arbitrary type,
and let $W(X)$ be the associated Coxeter group with distinguished
set of generating involutions $S(X)$.  Denote by $\H(X)$ the Hecke
algebra associated to $W(X)$.
(The reader is referred to \cite{{\bf 8}, \S7} for the basic theory of Hecke 
algebras arising from Coxeter systems.)  The $\A$-algebra $\H(X)$ has 
a basis consisting of elements $T_w$, with $w$ ranging over $W(X)$, 
that satisfy 
$$T_s T_w = \cases
T_{sw} & \text{ if } \ell(sw) > \ell(w),\cr
q T_{sw} + (q-1) T_w & \text{ if } \ell(sw) < \ell(w),\cr
\endcases$$ where $\ell$ is the length function on the Coxeter group
$W(X)$, $w \in W(X)$, and $s \in S(X)$.
The parameter $q$ is equal to $v^2$, where $\A = \zed[v, v^{-1}]$.

We define $J(X)$ to be the ideal of $\H(X)$ generated by 
all elements $$
\sum_{w \in \lan s_i, s_j \ran} T_w,
$$ where $(s_i, s_j)$ runs over all pairs of elements of $S(X)$
that correspond to adjacent nodes in the Coxeter graph.  
(If the nodes corresponding to $(s_i, s_j)$ are connected by a
bond of infinite strength, we omit the corresponding relation.)

\definition{Definition 1.1}
The generalized Temperley--Lieb algebra, $TL(X)$, is defined to be 
the quotient $\A$-algebra $\H(X)/J(X)$.
\enddefinition

The algebra $TL(X)$ has a basis which arises naturally from the
$T$-basis of $\H(X)$.  To define it, we introduce some standard 
combinatoric notions.

\definition{Definition 1.2}
Let $s_i$ and $s_j$ be elements of $S(X)$, the set of generating
involutions, such that $s_i$ and $s_j$ correspond to adjacent nodes in the
Coxeter graph.  Let $w_{ij}$ be the longest element in $\lan s_i, s_j
\ran$.  

We call an element $w \in W(X)$ {\it complex} if it can be written as
$x_1 w_{ij} x_2$, where $x_1, x_2 \in W(X)$ and $\ell(x_1 w_{ij} x_2)
= \ell(x_1) + \ell(w_{ij}) + \ell(x_2)$.

Denote by $W_c(X)$ the set of all elements of $W(X)$
which are not complex.

Let $t_w$ denote the image of the basis element $T_w \in \H(X)$ in
the quotient $TL(X)$.
\enddefinition

\proclaim{Theorem 1.3 (Graham)}
The set $\{ t_w : w \in W_c \}$ 
is an $\A$-basis for the algebra $TL(X)$.
\endproclaim

\demo{Proof}
This is \cite{{\bf 5}, Theorem 6.2}.
\qed\enddemo

We will call the basis of Theorem 1.3 the ``$t$-basis'' of the
algebra $TL(X)$.  It plays an important r\^ole in the sequel, since the
canonical basis will be defined in terms of it.

\proclaim{Lemma 1.4}
The algebra $TL(X)$ has a $\zed$-linear automorphism of order $2$
which sends $v$ to $v^{-1}$ and $t_w$ to $t_{w^{-1}}^{-1}$.
\endproclaim

\demo{Note}
The statement of the lemma makes sense: it is
well known that the elements $T_w \in \H(W)$ are invertible, from
which it follows that the elements $t_w$ are invertible also.
\enddemo

\demo{Proof}
It is known from \cite{{\bf 11}} that the Hecke algebra $\H(X)$ over
$\zed[v, v^{-1}]$ has a $\zed$-linear automorphism of order $2$ which
exchanges $v$ and
$v^{-1}$ and sends $T_w$ to $T_{w^{-1}}^{-1}$.  It is therefore
enough to show that this automorphism of $\H(X)$ fixes the ideal
$J(X)$.

Given a finite Coxeter group $W'$ with longest element $w_0$, it is
a standard result that, for any $w_1 \in W'$, there exists an element $w_2
\in W'$ such that $w_0 = w_1 w_2$ and $\ell(w_0) =
\ell(w_1) + \ell(w_2)$.  It follows from this fact that $$
\left(\sum_{w \in \lan s_i, s_j \ran} T_w\right) T_{w_{ij}}^{-1} =
\left(\sum_{w \in \lan s_i, s_j \ran} T_{w^{-1}}^{-1}\right)
,$$ so that the right hand side of the above equation lies in the
ideal $J(X)$.  We deduce that the automorphism of $\H(X)$ given above
fixes the ideal $J(X)$ setwise, and thus induces an automorphism of
$TL(X)$ which has the required properties.
\qed\enddemo

\proclaim{Lemma 1.5}
Let $w \in W$ (not necessarily in $W_c$).  Then $$
t_w = \sum_{{x \in W_c} \atop {x \leq w}}D_{x,w}t_x
,$$ where $D_{x, w} \in \zed[q]$.  Furthermore, $D_{w, w} = 1$ if $w
\in W_c$.
\endproclaim

\demo{Proof}
We proceed by induction on $\ell(w)$.
The proposition is trivial if $w \in W_c$, which covers the cases
where $\ell(w) \leq 1$, as well as the last assertion in the statement.

Now consider the case where $w = su$, $\ell(s) = 1$ and $\ell(u) =
\ell(w) - 1$.  Then, by induction, $$
t_w = t_s t_u = \sum_{{x \in W_c} \atop {x \leq u}}D_{x,u} t_s t_x
.$$ 

It follows from standard properties of $\H(X)$ that $$
t_s t_x = \cases
t_{sx} & \text{ if } \ell(sx) > \ell(x),\cr
q t_{sx} + (q-1) t_x & \text{ if } \ell(sx) < \ell(x).\cr
\endcases
$$  It is easy to see that $sx \leq w$ and $x \leq w$ for each basis
element $t_x$
appearing in the sum.  If $sx < w$ then, by induction, $t_{sx}$ is a
$\zed[q]$-linear combination of basis elements $t_z$ ($z \in W_c$) where 
$z \leq sx \leq w$.  

The only remaining case is when $sx = w$, $w \not\in W_c$.  This
forces $x = u$ and thus $u \in W_c$, so there is only one term in the
sum.  Since $w$ is
complex, it can be written as $x_1 w_{ij} x_2$ as in Definition 1.2.
We now use the fact that in $TL(X)$ we have $$
t_{w_{ij}} = - \sum_{y < w_{ij}} t_y
.$$  This enables $t_w = t_{x_1} t_{w_{ij}} t_{x_2}$ to be expressed
as a $\zed[q]$-linear
combination of elements $t_y$, where $y < w$.  By induction, these
elements $t_y$ are $\zed[q]$-linear combinations of basis elements
$t_z$ ($z \in W_c$) where $z \leq y \leq w$.

The lemma follows.
\qed\enddemo

\head 2. IC bases \endhead

We now recall the basic properties of IC bases from \cite{{\bf 2}, \S1}.

Let $\A := \zed[v, v^{-1}]$, where $v$ is an indeterminate, and let
$\A^- = \zed[v^{-1}]$ and $\A^+ = \zed[v]$.  Let $\,\bar{ }\,$ be the
involution on the ring $\A$ which satisfies $\bar{v} = v^{-1}$.

Let $M$ be a free $\A$-module with basis $\{m_i\}_{i \in I}$ and an
involutive $\zed$-linear map $\bar{ } : M \ra M$ such that
$\overline{am} = \bar{a} \bar{m}$ for any $m \in M$ and $a \in \A$.
For each $i \in I$, let $r_i$ be an integer.  

Let $\L$ be the free $\A^-$-submodule with basis $\{m'_i\}_{i \in I}$, 
where $m'_i := v^{r_i} m_i$.  Let $\pi : \L \ra \L/v^{-1}\L$ be the
canonical projection.

\definition{Definition 2.1}
If there exists a unique basis $\{ c_i \}_{i \in I}$ for $\L$
such that $\overline{c_i} = c_i$ and $\pi(c_i) = \pi(m'_i)$, then the
basis $\{c_i\}_{i \in I}$ is called an IC basis of $M$ with respect to
the triple $(\{m_i\}_{i \in I}, \ \bar{ } \ , \L)$.
\enddefinition

Note that the existence of an IC basis for $M$ depends only on the
original basis $\{m_i\}_{i \in I}$, the map $\,\bar{ }\,$ and the 
integers $r_i$.

\proclaim{Theorem 2.2 (Du)}
Let $(I, \leq)$ be a poset such that the sets $\{i : i \in I, i \leq
j\}$ are finite for all $j \in I$.  Suppose that $$
\overline{m'_j} = \sum_{{i \in I} \atop {i \leq j}} a_{ij} m'_i
$$ with $a_{ij} \in \A$ such that $a_{ii} = 1$ for all $i \in I$.
Then an IC basis of $M$ with respect
to $(\{m_i\}, \ \bar{ } \ , \L)$ exists.
\endproclaim

\demo{Proof}
This comes from \cite{{\bf 2}, Theorem 1.2, Remark 1.2.1 (1)}.
\qed\enddemo

The importance of this result for us is that the $t$-basis of the
algebra $TL(X)$ fits naturally into this setup, as we now explain.

Fix a Coxeter graph $X$.  Let $I = I(X)$ be the set $\{ w : w \in W_c
\}$.  We make $I$ into a poset $(I, \leq)$ by restricting the
Bruhat--Chevalley order on the Coxeter group $W = W(X)$ to the subset $I$.
(Standard properties of this order imply that the poset $(I,
\leq)$ has the finiteness property required by Theorem 2.2.)
We take $r_w = -\ell(w)$ and $m_w=t_w$, so that 
$m'_w = v^{-\ell(w)}t_w$, and we take $\,\bar{ }\,$ to be 
the automorphism of $TL(X)$ defined in Lemma 1.4.

Maintaining this notation, we have the following central result.

\proclaim{Theorem 2.3}
Let $X$ be an arbitrary Coxeter graph.
There exists an IC basis for the algebra $TL(X)$ with respect to 
$(\{m_i\}, \ \bar{ } \ , \L)$.
\endproclaim

\demo{Proof}
It is enough to show that the formula for $\overline{m'_w}$ 
(where $w\in W_c$) given
in Theorem 2.2 can be satisfied for suitable elements $a_{ij} \in \A$.
We proceed by induction on $\ell(w)$, the case $\ell(w) = 0$ being
trivial ($w = e$ and $\overline{m_e} = m_e$).

To deal with the inductive step, we suppose that $w = su$, 
where $\ell(s) = 1$ and $\ell(u) = \ell(w) - 1$.  
It is clear from the definition of $W_c$ that $u \in W_c$.

Using the fact that $m'_s = v^{-1} t_s$, it is easily verified that $$
\overline{m'_s} = m'_s - (v - v^{-1}) m'_e
,$$ where $m'_e$ is the identity in $TL(X)$.  Thus $$
\overline{m'_w} = \overline{m'_s} \overline{m'_u} = 
(m'_s - (v - v^{-1}) m'_e) \overline{m'_u}
.$$  By induction, $$
\overline{m'_u} = 
\sum_{{z \in W_c} \atop {z \leq u}} a_{zu} m'_z
.$$  For each $z$ appearing in the sum, we have $z \leq u \leq w$ and
also $sz \leq w$ by standard properties of the Bruhat--Chevalley order.
It therefore follows from Lemma 1.5 that $$
\overline{m'_w} = 
\sum_{{y \in W_c} \atop {y \leq w}} a_{yw} m'_y
$$ for suitable $a_{yw} \in \A$.

It remains to show that $a_{ww} = 1$.  Considering lengths and using 
Lemma 1.5, we find that the only
term in the expression for $\overline{m'_s} \overline{m'_u}$ which can
contribute to the coefficient of $m'_w$ arises from the product $m'_s
\times a_{uu} m'_u$.  It is clear that $m'_s m'_u = m'_w$, and
we have $a_{uu} = 1$ by induction, completing the proof.
\qed\enddemo

\demo{Remark 2.4}
\item{(1)}{The structure constants 
associated to the IC basis of $TL(X)$ lie in $\enn[v, v^{-1}]$
if $X$ is of type $A$, $D$ or $E_n$. (This will follow from Theorem 3.6.)
We do not know an example where this positivity property
fails for the IC basis, and the property 
certainly holds for some non-simply-laced types, such as $TL(H_n)$. \newline
In \cite{{\bf 7}, \S4.1}, the first author constructed a basis of
$TL(H_n)$ with structure constants in $\enn[v, v^{-1}]$ using
morphisms in the category of decorated tangles. (This category
was introduced in \cite{{\bf 6}}.)  It turns out that the basis of 
$TL(H_n)$ in \cite{{\bf 7}} is precisely the IC basis of $TL(H_n)$
arising from Theorem 2.3.  This means that the results of this
paper give an elementary characterisation of the basis in 
\cite{{\bf 7}}.  We do not supply the details here, but we hope to
do so in a forthcoming paper.
}
\item{(2)}{An interesting problem to consider is that of identifying the 
Coxeter graphs $X$
for which the IC basis of $TL(X)$ is equal to the image of the set 
$\Cal C$ of Kazhdan--Lusztig
basis elements $C'_w \in \H(X)$ indexed by $w \in W_c$.
Coxeter graphs $X$ of type $A$ have this property, by 
\cite{{\bf 4}, Theorem 3.8.2} and our Theorem 3.6.  
It seems likely that the set $\Cal C$ projects to the IC basis 
in type $D$, as well. \newline
This problem is closely related to the question of whether an
element $v^{-\ell(w)}T_w \in \H(X)$ ($w \notin W_c$) necessarily
lies in the lattice $\L$ after passing to $TL(X)$.  One can also 
express the problem in terms of a degree bound on the polynomials
$D_{x,w}$ of Lemma 1.5. \newline
A significant partial result would be to know that 
$v^{-\ell(w)}T_w \in \H(X)$ projects into $\L$ when $w = su$,
$s \in S(X)$, $u \in W_c$ and $w \notin W_c$.  This
would allow an inductive construction of the IC basis.
}
\enddemo

\head 3. The $ADE$ case \endhead

In \S3, we restrict ourselves to the case where the Coxeter graph is
of type $A$, $D$ or $E$.  However, we allow the graphs of type $E$ to
be of arbitrary rank.  This means that a graph of type $E_n$ ($n \geq
6$) may consist of a straight
line of nodes numbered $1, 2, \ldots, n-1$ together with a node
numbered $0$ which is connected only to node $3$.  (Contrast this with
the more familiar definition of type $E$ which is the same but
requires $n \leq 8$.)

It is known \cite{{\bf 5}, Theorem 7.1} 
that the algebras $TL(E_n)$ are finite dimensional for all values of
$n$.  In fact, for simply laced $X$, the algebra $TL(X)$ is finite 
dimensional if and only if $X$ is of type $A$, type $D$ or type $E_n$ 
for some $n$.  If $X$ satisfies these hypotheses, we say it is of type $ADE$.

\definition{Definition 3.1}
If $s \in S(X)$, we define $b_s$ to be the element $v^{-1} t_s +
v^{-1} t_e$.

If $w \in W_c$ and $w = s_1 s_2 \cdots s_r$ (reduced), we define $$
b_w := b_{s_1} b_{s_2} \cdots b_{s_r}.
$$
\enddefinition

This definition may appear to depend on the reduced expression chosen for
$w$, but in fact it does not, because $b_s$ and $b_{s'}$ commute whenever
$s$ and $s'$ commute, and every reduced expression for $w$ can 
be obtained from any other by a sequence of commutation moves 
(see \cite{{\bf 3}}).

It is well known (and follows easily from Theorem 1.3 and Lemma 1.5)
that the set $\{ b_w : w \in W_c \}$ is a basis for $TL(X)$.  We call
this the {\it monomial basis}.

We are particularly interested in the case where the Coxeter graph $X$
is connected and simply laced (all bonds have strength $2$ or $3$).

\proclaim{Lemma 3.2}
If $X$ is connected and simply laced then
the algebra $TL(X)$ is generated
as an associative, unital algebra by the elements $b_s$ (one for each
node of $X$) and defining relations $$\eqalign{
b_s^2 &= q_c b_s, \cr
b_s b_t &= b_t b_s \quad \text{ if $s$ and $t$ are not connected} \cr
b_s b_t b_s &= b_s \quad \text{ if $s$ and $t$ are connected} \cr
}$$ where $q_c := v + v^{-1}$.
\endproclaim

\demo{Proof}
This is a standard result from \cite{{\bf 3}}.
\qed\enddemo

Our aim in \S3 is to prove that in types $A$, $D$ and $E$, the monomial
basis is the IC basis.  The following is a key ingredient of the proof
that we will give.

\proclaim{Lemma 3.3}
Let $X$ be a Coxeter graph of type $ADE$.  Let $w \in W_c(X)$ and let 
$s_1s_2\cdots s_r$ be a reduced expression for $w$.
Then for any $1\leq i_1<i_2<\cdots <i_k\leq r$ ($k<r$), we have
$b_{s_{i_1}}b_{s_{i_2}}\cdots b_{s_{i_k}}=q_{c}^mb_x$ for some 
$x \in W_c(X)$, where $m\leq r-k-1$ and $q_c := v + v^{-1}$.
\endproclaim

\demo{Proof}
This result is a special case of \cite{{\bf 5}, Lemma 9.13}.  
Let $b$ be an arbitrary monomial in the generators $b_s$.  The relations of
$TL(X)$ ensure that $b$ is equal to $q_c^m b_x$ for some nonnegative
integer $m$ and some $x \in W_c$.  Then \cite{{\bf 5}, Lemma 9.13} shows 
that the removal of one generator from the monomial $b$ results
in an expression equal to $q_c^{m'} b_{x'}$, where $x' \in W_c$ and 
$m' \leq m + 1$.
Furthermore, if $b = b_w$ for some $w \in W_c$, we have $m' = m$.
This means that the
maximum exponent of $q_c$ which could occur after $r - k$ generators
have been removed from $b_w$, as in the statement, is $r - k - 1$.
\qed\enddemo

\demo{Remark 3.4}
It is possible to prove Lemma 3.3 without appealing to the results of
\cite{{\bf 5}} by using a combinatoric argument based on tom Dieck's
graphical calculus \cite{{\bf 1}}, and in fact the latter approach
establishes the lemma for
a slightly wider class of Coxeter systems, including type $\widetilde{E}_7$.
However, we do not pursue this here.
\enddemo

In the following lemma we use the standard notation $\e_x := (-1)^{\ell(x)}$.

\proclaim{Lemma 3.5}
Let $w \in W_c$.  Then $$
v^{-\ell(w)} t_w = \e_w \sum_{{x \in W_c} \atop {x \leq w}} \e_x
\Q_{x, w} b_x,
$$ where $\Q_{w, w} = 1$ and $\Q_{x, w} \in v^{-1} \A^-$ if $x < w$.
\endproclaim

\demo{Proof}  Fix $w \in W_c$ and fix a reduced expression $w = s_1
s_2 \cdots s_r$.  We have $$\eqalign{
v^{-\ell(w)} t_w
&= (v^{-1} t_{s_1})(v^{-1} t_{s_2}) \cdots (v^{-1} t_{s_r}) \cr
&= (b_{s_1} - v^{-1})(b_{s_2} - v^{-1}) \cdots (b_{s_r} - v^{-1}). \cr
}$$  By expanding the last product and using Lemma 3.2 and the
subexpression characterisation of the Bruhat--Chevalley order, we can
see that all the $b_x$ occurring in the sum satisfy $x \in W_c$ and $x
\leq w$.

More precisely, the product expands to a sum of terms $$
v^{k - \ell(w)} b_{s_{i_1}} b_{s_{i_2}} \cdots b_{s_{i_k}}.
$$  If $k < \ell(w)$, Lemma 3.3 shows that this is equal to $$
v^{k - \ell(w)}q_c^m b_x
,$$ where $x \in W_c$ and $m \leq \ell(w) - k - 1$.  It follows that
the coefficient of $b_x$, and hence $\Q_{x, w}$,
lies in $v^{-1} \A^-$ in this case.  The other
possibility is that $k = \ell(w)$, which produces the basis element
$b_w$ with coefficient $1$ and no other basis elements.
It follows that $\Q_{w, w} = 1$, as required.
\qed\enddemo

We are ready to prove the main result of \S3.

\proclaim{Theorem 3.6}
Let $X$ be a Coxeter graph of type $ADE$.  Then 
the IC basis for $TL(X)$ defined in Theorem 2.3
is precisely the monomial basis $\{b_w : w \in W_c\}$.
\endproclaim

\demo{Proof}
It is easily checked that the generators $b_s$ are fixed by the
involution $\,\bar{ }\,\,$, from which it follows that any monomial in 
these generators also has this property, so that $\overline{b_w} = b_w$.

Lemma 3.5 shows that the transfer matrix $(\e_w \e_x \Q_{x, w})$ 
from the monomial basis to the $m'_i$-basis is upper unitriangular
with respect to a total refinement of the partial order $\leq$, and
all the entries above the diagonal lie in $v^{-1} \A^-$.
Elementary linear algebra shows that the inverse of this matrix,
$(\P_{x, w})$, is also upper unitriangular with all the entries above
the diagonal in $v^{-1} \A^-$.  This shows that $$
b_w = \sum_{{x \in W_c} \atop {x \leq w}} \P_{x, w} (v^{-\ell(x)} t_x),
$$ and therefore that $\pi(b_w) = \pi(v^{-\ell(w)} t_w)$.

We have now shown that the monomial basis is the IC basis with respect to the
triple $(\{m_i\}_{i \in I}, \ \bar{ } \ , \L)$, as claimed.
\qed\enddemo

\demo{Remark 3.7}
\item{(1)}
{It is not true that the monomial basis of $TL(X)$ equals the IC basis 
for all simply laced $X$.  
For example, take $X$ to be the Coxeter graph
of type $\widetilde{A}_3$ consisting of four nodes connected by four
edges in the shape of a square.  If we number the nodes $1, 2, 3, 4$
around the square, the element $b_1 b_3 b_2 b_4 b_1 b_3$ is not an IC
basis element.  Similar remarks hold for $X$ of type
$\widetilde{A}_l$, where $l>3$ is odd.

If $X$ is non-simply-laced, the monomial and IC bases do not agree: any 
element in $W_c$ of the form $w = ss's$ (where $s, s' \in S(X)$) has the
property that $b_s b_{s'} b_s$ is not an IC basis element.
}
\item{(2)}
{Another problem to consider is that of determining the precise 
relationship between the elements
$\P_{x, w}$ in the proof of Theorem 3.6 and the Kazhdan--Lusztig 
polynomials $P_{x, w}$ of \cite{{\bf 11}}.
It is natural to be curious about this relationship, because our elements 
$\P_{x, w}$ play a r\^ole analogous to that of $v^{\ell(x) - \ell(w)}
P_{x, w}$ in \cite{{\bf 11}}.  Similarly, the elements $\Q_{x, w}$ are
analogous to the inverse Kazhdan--Lusztig polynomials.
}
\enddemo

\leftheadtext{}
\rightheadtext{}

\end